\documentclass[reqno,12pt,a4paper]{amsart}

\UseRawInputEncoding
\usepackage{amssymb}
\usepackage{bbm}
\usepackage{}

\usepackage{tikz}
\usetikzlibrary{matrix,arrows,calc,snakes,patterns,decorations.markings}

\usepackage{booktabs}
\usepackage[mathscr]{eucal}
\usepackage{graphics,epic}
\usepackage{amsfonts}
\usepackage{amscd}
\usepackage{latexsym}
\usepackage{amsmath,amssymb, amsthm, stmaryrd, bm, bbm}
\usepackage[all,2cell]{xy}
\usepackage{mathrsfs}
\usepackage{color}

\newtheorem*{theorem*}{Theorem}

\newtheorem*{conjecture*}{Conjecture}

\newtheorem*{question*}{Question}
\theoremstyle{remark}

\theoremstyle{definition}

\newtheorem*{notation*}{Notation}

\newcommand{\opname}[1]{\operatorname{\mathsf{#1}}}

\renewcommand{\mod}{\opname{mod}\nolimits}

\numberwithin{equation}{section}

\def\REF#1{\par\hangindent\parindent\indent\llap{#1\enspace}\ignorespaces}

\begin{document}
\begin{sloppypar}

\title[]{A construction of support $\tau$-tilting modules over $\tau$-tilting finite algebras}

\author{Yingying Zhang}

\address{Department of Mathematics, Huzhou University, Huzhou 313000, Zhejiang Province, P.R.China}

\email{yyzhang@zjhu.edu.cn}

\begin{abstract}
The notion of (semi)bricks, regarded as a generalization of (semi)simple modules, appeared in a paper of Ringel in 1976. In recent years, there have been several new developments motivated by links to $\tau$-tilting theory studied by Demonet-Iyama-Jasso and Asai. In this paper, we show how to glue semibricks along a recollement. As an application, we investigate the behavior of $\tau$-tilting finite under recollements of module categories of algebras. Moreover, we give some examples to show the construction of support $\tau$-tilting modules over $\tau$-tilting finite algebras by gluing semibricks via recollements.
\end{abstract}

\maketitle

\section*{0 Introduction}
A recollement of abelian categories is an exact sequence of abelian categories where both the inclusion functor and the quotient functor admit left and right adjoints. They first appeared in the construction of the category of perverse sheaves on a singular space by Beilinson, Bernstein and Deligne$^{[6]}$, arising from recollements of triangulated categories. Also recollements are quite an active subject widely studied by many authors (see [3, 6, 7, 10, 15, 17, 22] and so on). Beilinson, Bernstein and Deligne applied gluing techniques for simple modules and t-structures with respect to a recollement. Gluing co-t-structures was studied in [5]. Liu-Vit\'oria-Yang discussed gluing of silting objects via a recollement of bounded derived categories of finite dimensional algebras over a field with respect to the gluing of co-t-structures [15], recently, in [21] via gluing t-structures. Parra and Vit\'oria investigated the gluing of some basic properties of abelian categories (well-poweredness, Grothendieck's axiom AB3, AB4 and AB5, existence of a generator)$^{[18]}$.

In the representation theory of finite-dimensional algebras, the set of (semi)simple modules is fundamental. By Schur's Lemma, the endomorphism ring of a simple module is a division algebra and there exists no nonzero homomorphism between two nonisomorphic simple modules. As a generalization, a module is called a brick if its endomorphism ring is a division algebra. A set of isoclasses of pairwise Hom-orthogonal bricks is called a semibrick. It has long been studied in representation theory$^{[11][12][20]}$. Recently, Demonet-Iyama-Jasso gave the relation between bricks and $\tau$-rigid modules$^{[9]}$ and Asai investigated semibricks from the point of view of $\tau$-tilting theory$^{[1]}$.

Since simple modules can be glued via a recollement (see [6]), this leads to the natural question of how to glue semibricks. It turns out an answer to the question of gluing semibricks can be given more easily. The second goal of this paper is to show that the process of gluing semibricks allows to construct support $\tau$-tilting modules over $\tau$-tilting finite algebras in the middle category of a recollement out of support $\tau$-tilting modules in its outer categories. This paper is organized as follows.

In Section 1, we introduce some terminology and preliminary results needed throughout the paper.

In Section 2, we study gluing semibricks via a recollement. Precisely, we prove the following theorem.

{\bf Theorem 0.1}{\rm(Proposition 2.1 and Theorem 2.1)}\quad
Let $R(A, B, C)$ be a recollement (see Definition 1.4 for details). Then the functor $i_{*}$ maps semibricks in $\opname{mod}A$ to semibricks in $\opname{mod}B$. The functors $j_{!}, j_{*}$ and $j_{!*}$ map semibricks in $\opname{mod}C$ to semibricks in $\opname{mod}B$. There is an injection between sets of isomorphism classes of semibricks:\\
\centerline{$\{\text{semibricks in }\opname{mod}A\}\cup  \{\text{semibricks in } \opname{mod}C\}\rightarrow \{\text{semibricks in } \opname{mod}B\}$}\\ through a semibrick $\mathcal{S}_{L}\in \opname{mod}A$ and a semibrick $\mathcal{S}_{R}\in \opname{mod}C$ into $i_{*}(\mathcal{S}_{L})\sqcup j_{!*}(\mathcal{S}_{R})$.

Applying the above result, we observe that in a recollement of module categories over finite dimensional algebras, if the middle algebra is $\tau$-tilting finite, the other two algebras involved are also $\tau$-tilting finite. But the converse does not hold in general.

{\bf Theorem 0.2}{\rm(Theorem 2.2)}\quad
Let $R(A, B, C)$ be a recollement. If $B$ is $\tau$-tilting finite, then $A$ and $C$ are $\tau$-tilting finite.

In Section 3, we show the construction of support $\tau$-tilting modules over $\tau$-tilting finite algebras by gluing semibricks via recollements. Remark that for a $\tau$-tilting finite algebra, support $\tau$-tilting modules are in bijective correspondence with semibricks$^{[1]}$. Our result on the construction of support $\tau$-tilting modules are based on the following theorem.

{\bf Theorem 0.3}{\rm(Theorem 3.1)}\quad
Let $R(A, B, C)$ be a recollement and $B$ a $\tau$-tilting finite algebra. If $M_{A}$ and $M_{C}$ are respectively support $\tau$-tilting modules in $\opname{mod}A$ and $\opname{mod}C$, with the corresponding semibricks $\mathcal{S}_{A}$ and $\mathcal{S}_{C}$, then there exists a unique support $\tau$-tilting $B$-module $M_{B}$ which is associated with the induced semibrick $i_{*}(\mathcal{S}_{A})\sqcup j_{!*}(\mathcal{S}_{C})$.

Finally, we give examples to illustrate the process.

{\bf Notation}\quad
Let $K$ be a field and $A$ a finite-dimensional $K$-algebra. We denote by $\opname{mod}A$(resp, $\opname{proj}A$) the category of finitely generated(resp, finitely generated projective) right $A$-modules and by $\tau$ the Auslander-Reiten translation of $A$. The composition of maps or functors $f:X\longrightarrow Y$ and $g:Y\longrightarrow Z$ is denoted by $gf$. For $M\in \mod A$, we denote by $\opname{ind}M$ the set of isoclasses of indecomposable direct summands of $M$. We use $``\sqcup"$ to denote disjoint union.

\section*{1 Preliminaries}

In this section, we collect some basic materials that will be used later. We begin with the definition of semibricks.

\vspace{0.2cm}
{\bf 1.1. Semibricks}

\vspace{0.2cm}
{\bf Definition 1.1}\quad
(1) A module $S\in \opname{mod}A$ is called a $brick$ if  $\opname{End}_{A}(S)$ is a division $K$-algebra (i.e., the non-trivial endomorphisms are invertible). We write $\opname{brick}A$ for the set of isoclasses of bricks in $\opname{mod}A$.\\
(2) A subset $\mathcal{S}\subset \opname{mod}A$ of isoclasses of bricks is called a $semibrick$ if $\opname{Hom}_{A}(S_{1}, S_{2})=0$ for any $S_{1}\neq S_{2}\in\mathcal{S}$. We write $\opname{sbrick}A$ for the set of semibricks in $\opname{mod}A$.

A typical example is that a simple module is a brick and a set of isoclasses of simple modules is a semibrick by Schur's Lemma. Moreover, preprojective modules and preinjective modules over a finite-dimensional algebra are also bricks in the module categories(see [4, \uppercase\expandafter{\romannumeral8}.2.7. Lemma]).

\vspace{0.2cm}

{\bf 1.2. $\tau$-tilting theory}

\vspace{0.2cm}

Recall the definition of support $\tau$-tilting modules from [2].

{\bf Definition 1.2}\quad
Let ($X$,$P$) be a pair with $X \in \mod A$ and $P \in \opname{proj}A$.\\
(1) We call $X$ in $\mod A$ $\tau$-$rigid$ if  ${\opname{Hom}}_{A}(X, \tau X)$=0. We call ($X$,$P$) a $\tau$-$rigid$  $pair$ if $X$ is $\tau$-rigid and ${\opname{Hom}}_{A}(P, X)$=0.\\
(2) We call $X$ in $\mod A$ $\tau$-$tilting$ if $X$ is $\tau$-rigid and $ |X| = |A|$, where $ |X|$ denotes the number of nonisomorphic indecomposable direct summands of $X$.\\
(3) We call $X$ in $\mod A$ $support$ $\tau$-$tilting$ if there exists an idempotent $e$ of $A$ such that $X$ is a $\tau$-tilting ($A/\langle e\rangle$)-module. We call ($X$,$P$) a $support$ $\tau$-$tilting$ $ pair$ if ($X$,$P$) is $\tau$-rigid and $|X|+|P|=|A|$.

We say that ($X$,$P$) is {\it basic} if $X$ and $P$ are basic. Moreover, for a support $\tau$-tilting pair ($X$,$P$), $X$ determines $P$ uniquely up to isomorphism. We denote by $\opname{s\tau-tilt\,}A$ the set of isomorphism classes of basic support $\tau$-tilting $A$-modules.

{\bf Definition 1.3}$^{[9]}$\quad
A finite-dimensional algebra $A$ is called $\tau$-$tilting$ $finite$ if there are only finitely many isomorphism classes of basic $\tau$-tilting $A$-modules.

In fact, by [9], $A$ is $\tau$-tilting finite if and only if there exist only finitely many isomorphism classes of indecomposable $\tau$-rigid $A$-modules, if and only if there exist only finitely many isomorphism classes of basic support $\tau$-tilting $A$-modules, if and only if there are only finitely many isomorphism classes of bricks in $\opname{mod}A$.

Note that Asai's result gave a bijection between support $\tau$-tilting modules and semibricks satisfying left finiteness condition, which is automatic for $\tau$-tilting finite algebras. We will not explain the finiteness condition, since we only consider the $\tau$-tilting finite algebras in section 3.

{\bf Theorem 1.1}\rm ([1,Theorem 1.3(2)])\quad
Let $A$ be a $\tau$-tilting finite algebra and $B:=\opname{End}_{A}(M)$, where $M\in \opname{mod}A$. There exists a bijection
\begin{center}
$\opname{s\tau-tilt\,}A\longrightarrow \opname{sbrick}A$
\end{center}
given by $M\mapsto  \opname{ind}(M/\opname{rad}_{B}M) $.

\vspace{0.2cm}

{\bf 1.3. Recollements}

\vspace{0.2cm}

For the convenience, we will recall the definition of recollements of abelian categories, see for instance [6,10,13].

{\bf Definition 1.4}\quad Let $A, B$ and $C$ be finite-dimensional algebras. Then a recollement of $\opname{mod}B$ relative to $\opname{mod}A$ and $\opname{mod}C$, diagrammatically expressed by
$$\xymatrix@!C=2pc{
\opname{mod}A\; \ar@{>->}[rr]|{i_{*}} && \opname{mod}B \ar@<-4.0mm>@{->>}[ll]_{i^{*}} \ar@{->>}[rr]|{j^{*}} \ar@{->>}@<4.0mm>[ll]^{i^{!}}&& \opname{mod}C \ar@{>->}@<-4.0mm>[ll]_{j_{!}} \ar@{>->}@<4.0mm>[ll]^{j_{*}}
}$$
which satisfies the following three conditions:\\
(1) ($i^{*}, i_{*}$), ($i_{*}, i^{!}$), ($j_{!}, j^{*}$) and ($j^{*}, j_{*}$) are adjoint pairs;\\
(2) $i_{*}, j_{!}$ and $j_{*}$ are fully faithful functors;\\
(3) ${\rm Im}i_{*}={\rm Ker}j^{*}$.

{\bf Remark}\quad
(1) From Definition 1.4(1), it follows that $i_{*}$ and $j^{*}$ are both right adjoint functors and left adjoint functors, therefore they are exact functors of abelian categories.

(2) By the definition of recollements, it follows that $i^{*}i_{*}\cong id, i^{!}i_{*}\cong id, j^{*}j_{!}\cong id$ and $j^{*}j_{*}\cong id$. Also $i^{*}j_{!}=0, i^{!}j_{*}=0$.

(3) Throughout this paper, we denote by $R(A, B, C)$ a recollement of $\opname{mod}B$ relative to $\opname{mod}A$ and $\opname{mod}C$ as above.

Associated to a recollement there is a seventh funtor $j_{!*}:=\opname{Im}(j_{!}\rightarrow j_{*}):\opname{mod}C\rightarrow \opname{mod}B$ called the intermediate extension functor. This functor plays an important role in gluing simple modules in [6]. Moreover, Crawley-Boevey-Sauter gave geometric applications$^{[8]}$; Keller and Scherotzke discussed stable and costable objects associated to a recollement (for the case of intermediate Kan extension)$^{[15]}$. The following proposition summarize results in [6,10].

{\bf Proposition 1.1}\quad
(1) In any recollement situation, we have $i^{*}j_{!*}=0, i^{!}j_{!*}=0$.\\
(2) $j^{*}j_{!*}\cong id$ and the functor $j_{!*}$ is fully faithful.\\
(3)The functor $j_{!*}$ sends simples in $\opname{mod}C$ to simples in $\opname{mod}B$. There is a bijection between sets of isomorphism classes of simples:
$$\{\text{simples in }\opname{mod}A\}\sqcup \{\text{simples in } \opname{mod}C\}\rightarrow \{\text{simples in } \opname{mod}B\}$$ given by mapping a simple $M_{L}\in \opname{mod}A$ to $i_{*}(M_{L})$ and a simple $M_{R}\in \opname{mod}C$ to $j_{!*}(M_{R})$.

\section*{2 Gluing semibricks}

In this section, we will discuss how to glue semibricks of the left category and the right category into a semibrick in the middle term with respect to a recollement.

{\bf Lemma 2.1}\quad
If $F:\opname{mod}A\rightarrow \opname{mod}B$ is a fully faithful functor, then we have $F(\opname{brick}A)\subseteq \opname{brick}B$ and $F(\opname{sbrick}A)\subseteq \opname{sbrick}B$.

{\bf Proof}\quad
(1) For any $S\in \opname{brick}A, \opname{End}_{A}(S)$ is a division ring. If we want to know $F(S)\in \opname{brick}B$, We only have to prove that $\opname{End}_{B}(F(S))$ is also a division ring. Let $0\neq h\in\opname{End}_{B}(F(S))$. Since $F$ is fully faithful, there exists $0\neq g\in\opname{End}_{A}(S)$ such that $F(g)=h$. It follows that there exists $g^{-1}\in\opname{End}_{A}(S)$ such that $g^{-1}g=gg^{-1}=id_{S}$, for the reason that $\opname{End}_{A}(S)$ is a division ring. Thus there exists $F(g^{-1})\in \opname{End}_{B}(F(S))$ such that $F(g^{-1})h=hF(g^{-1})=id_{F(S)}$.

(2) Let $\mathcal{S}\in \opname{sbrick}A$. By (1), $F(\mathcal{S})\subset \opname{brick}B$. Take $S_{1}, S_{2}\in \mathcal {S}$ such that $F(S_{1})\neq F(S_{2})$, then we have $S_{1}\neq S_{2}$. Since $F$ is fully faithful, $F(S_{1})\neq F(S_{2})$ satisfy $\opname{Hom}_{B}(F(S_{1}),F(S_{2}))\cong \opname{Hom}_{A}(S_{1},S_{2})= 0$. Therefore $F(\mathcal{S})\in\opname{sbrick}B$.

By Proposition 1.1(2) and Lemma 2.1, immediately we have

{\bf Proposition 2.1}\quad
Let $R(A, B, C)$ be a recollement.\\
(1)$i_{*}(\opname{brick}A)\subseteq\opname{brick}B$ and $i_{*}(\opname{sbrick}A)\subseteq\opname{sbrick}B$;\\
(2)$j_{!}(\opname{brick}C)\subseteq\opname{brick}B$ and $j_{!}(\opname{sbrick}C)\subseteq\opname{sbrick}B$;\\
(3)$j_{*}(\opname{brick}C)\subseteq\opname{brick}B$ and $j_{*}(\opname{sbrick}C)\subseteq\opname{sbrick}B$;\\
(4)$j_{!*}(\opname{brick}C)\subseteq\opname{brick}B$ and $j_{!*}(\opname{sbrick}C)\subseteq\opname{sbrick}B$.

Motivated by gluing simples in [6] (see Proposition 1.1(3)), here we also give a construction of semibricks from the left and right side into the middle.

{\bf Theorem 2.1}\quad
Let $R(A, B, C)$ be a recollement. We have
$i_{*}(\opname{sbrick}A)\sqcup j_{!*}(\opname{sbrick}C)\subseteq \opname{sbrick}B$.

{\bf Proof}\quad
Let $\mathcal{S}_{L}\in \opname{sbrick}A$ and $\mathcal{S}_{R}\in \opname{sbrick}C$. We claim that $i_{*}(\mathcal{S}_{L})\sqcup j_{!*}(\mathcal{S}_{R})\in\opname{sbrick}B$. By Proposition 2.1, $i_{*}(\mathcal{S}_{L})\sqcup j_{!*}(\mathcal{S}_{R})\subset\opname{brick}B$.
Let $S\neq S'\in i_{*}(\mathcal{S}_{L})\sqcup j_{!*}(\mathcal{S}_{R})$.\\
Case 1: Assume that $S\neq S'\in i_{*}(\mathcal{S}_{L})$.\\ Take $S=i_{*}(S_{l})$ and $S'=i_{*}(S_{l}')$ where $S_{l}\neq S_{l}'\in\mathcal{S}_{L}$.
Since $i_{*}$ is fully faithful, we have $\opname{Hom}_{B}(i_{*}(S_{l}), i_{*}(S_{l}'))\cong \opname{Hom}_{A}(S_{l},S_{l}')$. By $\mathcal{S}_{L}\in \opname{sbrick}A$, it follows that $\opname{Hom}_{B}(S,S')=0$.\\
Case 2: Assume that $S\neq S'\in j_{!*}(\mathcal{S}_{R})$. \\
Take $S=j_{!*}(S_{r})$ and $S'=j_{!*}(S_{r}')$ where $S_{r}\neq S_{r}'\in\mathcal{S}_{R}$.
Since $j_{!*}$ is fully faithful, we have $\opname{Hom}_{B}(j_{!*}(S_{r}), j_{!*}(S_{r}'))\cong \opname{Hom}_{C}(S_{r},S_{r}')$. By $\mathcal{S}_{R}\in \opname{sbrick}C$, it follows that $\opname{Hom}_{B}(S,S')=0$.\\
Case 3: Assume that $S\in i_{*}(\mathcal{S}_{L}), S'\in j_{!*}(\mathcal{S}_{R})$ or $S\in j_{!*}(\mathcal{S}_{R}), S'\in i_{*}(\mathcal{S}_{L})$.
By Definition 1.4(1) and Proposition 1.1(1), it follows that $$\opname{Hom}_{B}(i_{*}(\mathcal{S}_{L}), j_{!*}(\mathcal{S}_{R}))\cong \opname{Hom}_{A}(\mathcal{S}_{L}, i^{!}j_{!*}(\mathcal{S}_{R}))=0,$$ $$\opname{Hom}_{B}(j_{!*}(\mathcal{S}_{R}), i_{*}(\mathcal{S}_{L}))\cong \opname{Hom}_{A}(i^{*}j_{!*}(\mathcal{S}_{R}),\mathcal{S}_{L})=0.$$

The following definition settles what we will mean by gluing semibricks.

{\bf Definition 2.1}\quad
Let $R(A, B, C)$ be a recollement. We say that a semibrick $\mathcal{S}\in\mod B$ is glued from $\mathcal{S}_{L}\in\mod A$ and $\mathcal{S}_{R}\in\mod C$ with respect to $R(A, B, C)$ if $\mathcal{S}$ is obtained by the construction of Theorem 2.1.

By Proposition 2.1 and Theorem 2.1, we find ways to get some (not all) semibricks in the middle cateogry, but we don't know how to obtain the others. For example, the module category of Kronecker algebra is a recollement of two copies of the category of vector spaces, the construction in Theorem 2.1 only yields three semibricks but there are infinite semibricks in the module categories of Kronecker algebra.

Considering that $i^{*}j_{!}$ and $i^{!}j_{*}$ are both zero but $i^{!}j_{!}$ and $i^{*}j_{*}$ may not be zero in a recollement, $i_{*}(\opname{sbrick}A)\sqcup j_{!}(\opname{sbrick}C)$ and $i_{*}(\opname{sbrick}A)\sqcup j_{*}(\opname{sbrick}C)$ will not usually be semibricks in the middle category. Thanks to the result of [10, Proposition 8.8], we can still find the special case when $i^{!}j_{!}=0$ or $i^{*}j_{*}=0$. The same as the proof of Theorem 2.1, immediately we have the following corollary.

{\bf Corollary 2.1}\quad For a recollement, if $i^{*}$ is exact, we have $i_{*}(\opname{sbrick}A)\sqcup j_{!}(\opname{sbrick}C)\subseteq \opname{sbrick}B$. Dually, if $i^{!}$ is exact, we have $i_{*}(\opname{sbrick}A)\sqcup j_{*}(\opname{sbrick}C)\subseteq \opname{sbrick}B$.

As an application of Theorem 2.1, we give the behavior of $\tau$-tilting finite under recollements. According to [19], a recollement whose terms are module categories of finite-dimensional algebras is equivalent to one induced by an idempotent element. We should remark that the following result is a well-known result (see [16, Corollaries 2.3 and 2.4]).

{\bf Theorem 2.2}\quad
Let $R(A, B, C)$ be a recollement. If $B$ is $\tau$-tilting finite, then $A$ and $C$ are $\tau$-tilting finite.

{\bf Proof}\quad
Assume that $B$ is $\tau$-tilting finite. By Theorem 1.1, the set $\opname{sbrick}B$ is finite. According to Theorem 2.1, the sets $i_{*}(\opname{sbrick}A)$ and $j_{!*}(\opname{sbrick}C)$ are both finite. By Definition 1.4 and Proposition 1.1(2), $i_{*}$ and $j_{!*}$ are fully faithful, thus the sets $\opname{sbrick}A$ and $\opname{sbrick}C$ are finite. Clearly, the sets $\opname{brick}A$ and $\opname{brick}C$ are also finite. Applying [9, Theorem 1.4], we have finished to prove that $A$ and $C$ are $\tau$-tilting finite.

\section*{3 A construction of support $\tau$-tilting modules}
\medskip
Throughout this section, $R(A,B,C)$ is a recollement of module categories and $B$ is a $\tau$-tilting finite algebra. Since semibricks can be glued via a recollement, the natural question is the following:

{\bf Question}\quad
Given a recollement of module categories, support $\tau$-tilting modules $M_{A}$ in $\mod A$ and $M_{C}$ in $\mod C$, is it possible to construct a support $\tau$-tilting module in $\mod B$ corresponding to the glued semibrick?

In the following, we give a positive answer to this question. The idea of proof is to use the bijection between support $\tau$-tilting modules and semibricks over $\tau$-tilting finite algebras.

{\bf Theorem 3.1}\quad
If $M_{A}$ and $M_{C}$ are respectively support $\tau$-tilting modules in $\opname{mod}A$ and $\opname{mod}C$, with the corresponding semibricks $\mathcal{S}_{A}$ and $\mathcal{S}_{C}$, then there exists a unique support $\tau$-tilting $B$-module $M_{B}$ which is associated with the induced semibrick $i_{*}(\mathcal{S}_{A})\sqcup j_{!*}(\mathcal{S}_{C})$.

{\bf Proof}\quad
Assume that $M_{A}\in \opname{s\tau-tilt\,}A$ and $M_{C}\in \opname{s\tau-tilt\,}C$. Since $B$ is $\tau$-tilting finite, by Theorem 2.2 it follows that $A$ and $C$ are $\tau$-tilting finite. By Theorem 1.1, there exist the corresponding semibricks $\mathcal{S}_{A}\in\opname{sbrick}A$ and $\mathcal{S}_{C}\in\opname{sbrick}C$. Applying Theorem 2.1, we get the glued semibrick $i_{*}(\mathcal{S}_{A})\sqcup j_{!*}(\mathcal{S}_{C})\in\opname{sbrick}B$. Thus by Theorem 1.1 again, there exists a unique support $\tau$-tilting $B$-module $M_{B}$ which is associated with the induced semibrick $i_{*}(\mathcal{S}_{A})\sqcup j_{!*}(\mathcal{S}_{C})$.

Although we are convinced of the unique existence of the gluing of support $\tau$-tilting modules by gluing semibricks over $\tau$-tilting finite algebras, from the examples we can see that it is impossible to describe the construction of support $\tau$-tilting modules through a formula along a recollement(using only the modules one wants to glue and the functors of the recollement). From the left side, it is not even the functor $i_{*}$ (see Example 3.2). For the convenience of the reader, we give some examples to illustrate the process.

{\bf Example 3.1}\quad
Let $A$ be the path algebra over a field $K$ of the quiver $1\rightarrow 2\rightarrow 3$, of type $A_{3}$. If $e$ is the idempotent $e_{1}+ e_{2}$, then as a right $A$-module $A/\langle e\rangle $ is isomorphic to $S_{3}$ and $eAe$ is the path algebra of the quiver $1\rightarrow 2$. In this case, there is a recollement as follows:
$$\xymatrix@!C=2pc{
{\rm mod\,}(A/\langle e\rangle )\ar@{>->}[rr]|{i_{*}} && {\rm mod\,}A \ar@<-4.0mm>@{->>}[ll]_{i^{*}} \ar@{->>}[rr]|{j^{*}} \ar@{->>}@<4.0mm>[ll]^{i^{!}}&& {\rm mod\,}(eA e)\ar@{>->}@<-4.0mm>[ll]_{j_{!}} \ar@{>->}@<4.0mm>[ll]^{j_{*}}
}$$
where $i^{*}=-\otimes_{A}A/\langle e\rangle $, $j_{!}=-\otimes_{eA e}e A$, $i^{!}={\rm Hom}_{A}(A/\langle e\rangle , -)$, $i_{*}=-\otimes_{A/\langle e\rangle }A/\langle e \rangle $, $j^{*}=-\otimes_{A}Ae$, $j_{*}={\rm Hom}_{eAe}(Ae, -)$.

\begin{table}[htbp]\centering
 \begin{tabular}{lcl}
  \toprule
 $\opname{s\tau-tilt\,}(A/\langle e \rangle ) $\quad\quad & $\opname{s\tau-tilt\,}A$ \quad\quad\quad & $\opname{s\tau-tilt\,}(eAe)$ \\
  \midrule
  \quad\quad $\begin{smallmatrix}\color{red}{3}\end{smallmatrix}$ &
  $\begin{smallmatrix}\color{red}{3}\end{smallmatrix}\begin{smallmatrix} \color{red}{2}\\3 \end{smallmatrix}\begin{smallmatrix}\color{red} {1}\\2\\3 \end{smallmatrix}$
  \quad\quad\quad & \quad\quad
  $\begin{smallmatrix}{\color{red}{1}}\\2\end{smallmatrix}
  \begin{smallmatrix}\color{red}{2}\end{smallmatrix}$ \\\\
  \quad\quad $\begin{smallmatrix}\color{red}{3}\end{smallmatrix}$&
  $\begin{smallmatrix}\color{red}{3}\end{smallmatrix}\begin{smallmatrix} \color{red}{1}\\ \color{red}{2}\\3 \end{smallmatrix}\begin{smallmatrix}1 \end{smallmatrix}$\quad\quad\quad &
  \quad\quad
  $\begin{smallmatrix}\color{red}{1}\\ \color{red}{2}\end{smallmatrix}
  \begin{smallmatrix}1\end{smallmatrix}$\\\\
  \quad\quad $\begin{smallmatrix}\color{red}{3}\end{smallmatrix}$ &
  $\begin{smallmatrix}\color{red}{3}\end{smallmatrix}\begin{smallmatrix} \color{red}{2}\\3 \end{smallmatrix}$
  \quad\quad\quad & \quad\quad
  $\begin{smallmatrix}\color{red}{2}\end{smallmatrix}$ \\\\
 \quad\quad $\begin{smallmatrix}\color{red}{3}\end{smallmatrix}$ &
  $\begin{smallmatrix}\color{red}{3}\end{smallmatrix}\begin{smallmatrix} \color{red}{1} \end{smallmatrix}$
  \quad\quad\quad & \quad\quad
  $\begin{smallmatrix}\color{red}{1}\end{smallmatrix}$ \\\\
  \quad\quad $\begin{smallmatrix}\color{red}{3}\end{smallmatrix}$ &
  $\begin{smallmatrix}\color{red}{3}\end{smallmatrix}$
  \quad\quad\quad & \quad\quad
  $\begin{smallmatrix}0\end{smallmatrix}$ \\\\
  \quad\quad $\begin{smallmatrix}0\end{smallmatrix}$ &
  $\begin{smallmatrix}{\color{red}{1}}\\2\end{smallmatrix}
  \begin{smallmatrix}\color{red}{2}\end{smallmatrix}$
  \quad\quad\quad & \quad\quad
  $\begin{smallmatrix}{\color{red}{1}}\\2\end{smallmatrix}
  \begin{smallmatrix}\color{red}{2}\end{smallmatrix}$ \\\\
  \quad\quad $\begin{smallmatrix}0\end{smallmatrix}$ &
  $\begin{smallmatrix}{\color{red}{1}}\\\color{red}{2}\end{smallmatrix}
  \begin{smallmatrix}1\end{smallmatrix}$
  \quad\quad\quad & \quad\quad
  $\begin{smallmatrix}{\color{red}{1}}\\\color{red}{2}\end{smallmatrix}
  \begin{smallmatrix}1\end{smallmatrix}$ \\\\
  \quad\quad $\begin{smallmatrix}0\end{smallmatrix}$ &
  $\begin{smallmatrix}\color{red}{2}\end{smallmatrix}$
  \quad\quad\quad & \quad\quad
  $\begin{smallmatrix}\color{red}{2}\end{smallmatrix}$ \\\\
  \quad\quad $\begin{smallmatrix}0\end{smallmatrix}$ &
  $\begin{smallmatrix}\color{red}{1}\end{smallmatrix}$
  \quad\quad\quad & \quad\quad
  $\begin{smallmatrix}\color{red}{1}\end{smallmatrix}$ \\\\
  \quad\quad $\begin{smallmatrix}0\end{smallmatrix}$ &
  $\begin{smallmatrix}0\end{smallmatrix}$
  \quad\quad\quad & \quad\quad
  $\begin{smallmatrix}0\end{smallmatrix}$ \\\\
  \bottomrule
 \end{tabular}
\end{table}

In the above table, we give a complete list of support $\tau$-tilting module $M\in\mod A$ constructed from $M_{L}\in\mod (A/\langle e\rangle )$ and $M_{R}\in\mod (eAe)$, also the bricks in the corresponding semibrick for each support $\tau$-tilting are red in color.

The number of support $\tau$-tilting $A$ modules is fourteen (see [1, Example 1.18]). Note that there are ten support $\tau$-tilting $A$ modules in the middle line of the table and we cannot obtain the others by the construction.

{\bf Example 3.2}\quad
Let $A$ be the preprojective algebra of type $A_{3}$ which is given by the following quiver and relation $aa'=0, b'b=0, bb'=a'a$.
$$\xymatrix{
1 \ar@/^/[r]^{a}&
2 \ar@/^/[l]^{a'}\ar@/^/[r]^{b}&
3 \ar@/^/[l]^{b'}&
}$$
Let $e=e_{1}+ e_{3}$. Then as a right $A$-module $A/\langle e\rangle $ is isomorphic to $S_{2}$ and $eAe$ is the preprojective algebra of type $A_{2}$. Then there is a recollement $R(A/\langle e\rangle, A, eAe)$ induced by the idempotent $e$.

\begin{table}[htbp]\centering
 \begin{tabular}{lcl}
  \toprule
 $\opname{s\tau-tilt\,}(A/\langle e \rangle ) $\quad\quad & $\opname{s\tau-tilt\,}A$ \quad\quad\quad & $\opname{s\tau-tilt\,}(eAe)$ \\
  \midrule
  \quad\quad $\begin{smallmatrix}\color{red}{2}\end{smallmatrix}$ &
  $\begin{smallmatrix}\color{red}{1}\\ 2\\3\end{smallmatrix}\begin{smallmatrix} \color{red}{2}\\13\\2 \end{smallmatrix}\begin{smallmatrix}\color{red} {3}\\2\\1 \end{smallmatrix}$
  \quad\quad\quad & \quad\quad
  $\begin{smallmatrix}{\color{red}{1}}\\3\end{smallmatrix}
  \begin{smallmatrix}\color{red}{3}\\1\end{smallmatrix}$ \\\\
  \quad\quad $\begin{smallmatrix}\color{red}{2}\end{smallmatrix}$&
  $\begin{smallmatrix}\color{red}{2}\\3\end{smallmatrix}\begin{smallmatrix} 3\\2 \end{smallmatrix}\begin{smallmatrix}\color{red}{3}\\\color{red}{2}\\ \color{red}{1} \end{smallmatrix}$\quad\quad\quad &
  \quad\quad
  $\begin{smallmatrix}3\end{smallmatrix}
  \begin{smallmatrix}\color{red}{3}\\ \color{red}{1}\end{smallmatrix}$\\\\
  \quad\quad $\begin{smallmatrix}\color{red}{2}\end{smallmatrix}$ &
  $\begin{smallmatrix}\color{red}{1}\\ \color{red}{2}\\ \color{red}{3}\end{smallmatrix}\begin{smallmatrix} 1\\2 \end{smallmatrix}\begin{smallmatrix}\color{red}{2}\\ 1\end{smallmatrix}$
  \quad\quad\quad & \quad\quad
  $\begin{smallmatrix}\color{red}{1}\\ \color{red}{3}\end{smallmatrix}\begin{smallmatrix}1\end{smallmatrix}$ \\\\
 \quad\quad $\begin{smallmatrix}\color{red}{2}\end{smallmatrix}$ &
  $\begin{smallmatrix}\color{red}{2}\\3\end{smallmatrix}\begin{smallmatrix} \color{red}{3} \\2 \end{smallmatrix}$
  \quad\quad\quad & \quad\quad
  $\begin{smallmatrix}\color{red}{3}\end{smallmatrix}$ \\\\
  \quad\quad $\begin{smallmatrix}\color{red}{2}\end{smallmatrix}$ &
  $\begin{smallmatrix}\color{red}{1}\\2\end{smallmatrix}\begin{smallmatrix} \color{red}{2} \\1 \end{smallmatrix}$
  \quad\quad\quad & \quad\quad
  $\begin{smallmatrix}\color{red}{1}\end{smallmatrix}$ \\\\
  \quad\quad $\begin{smallmatrix}\color{red}{2}\end{smallmatrix}$ &
  $\begin{smallmatrix}\color{red}{2}\end{smallmatrix}$
  \quad\quad\quad & \quad\quad
  $\begin{smallmatrix}0\end{smallmatrix}$ \\\\
  \quad\quad $\begin{smallmatrix}0\end{smallmatrix}$ &
  $\begin{smallmatrix}{\color{red}{3}}\end{smallmatrix}
  \begin{smallmatrix}\color{red}{1}\end{smallmatrix}$
  \quad\quad\quad & \quad\quad
  $\begin{smallmatrix}{\color{red}{1}}\\3\end{smallmatrix}
  \begin{smallmatrix}\color{red}{3}\\1\end{smallmatrix}$ \\\\
  \quad\quad $\begin{smallmatrix}0\end{smallmatrix}$ &
  $\begin{smallmatrix}3\end{smallmatrix}\begin{smallmatrix}3\\2\end{smallmatrix}
  \begin{smallmatrix}\color{red}3\\ \color{red}2\\{\color{red}{1}}\end{smallmatrix}$
  \quad\quad\quad & \quad\quad
  $\begin{smallmatrix}3\end{smallmatrix}
  \begin{smallmatrix}\color{red}3\\ \color{red}1\end{smallmatrix}$ \\\\
  \quad\quad $\begin{smallmatrix}0\end{smallmatrix}$ &
  $\begin{smallmatrix}\color{red}{1}\\ \color{red}2\\ \color{red}3\end{smallmatrix}\begin{smallmatrix}1\\2\end{smallmatrix}
  \begin{smallmatrix}1\end{smallmatrix}$
  \quad\quad\quad & \quad\quad
  $\begin{smallmatrix}\color{red}{1}\\ \color{red}3\end{smallmatrix}\begin{smallmatrix}1\end{smallmatrix}$ \\\\
  \quad\quad $\begin{smallmatrix}0\end{smallmatrix}$ &
  $\begin{smallmatrix}\color{red}{3}\end{smallmatrix}$
  \quad\quad\quad & \quad\quad
  $\begin{smallmatrix}\color{red}{3}\end{smallmatrix}$ \\\\
  \quad\quad $\begin{smallmatrix}0\end{smallmatrix}$ &
  $\begin{smallmatrix}\color{red}1\end{smallmatrix}$
  \quad\quad\quad & \quad\quad
  $\begin{smallmatrix}\color{red}1\end{smallmatrix}$ \\\\
  \quad\quad $\begin{smallmatrix}0\end{smallmatrix}$ &
  $\begin{smallmatrix}0\end{smallmatrix}$
  \quad\quad\quad & \quad\quad
  $\begin{smallmatrix}0\end{smallmatrix}$ \\\\
  \bottomrule
 \end{tabular}
\end{table}
In the above table, we give a complete list of support $\tau$-tilting module $M\in\mod A$ constructed from $M_{L}\in\mod (A/\langle e\rangle )$ and $M_{R}\in\mod (eAe)$, also the bricks in the corresponding semibrick for each support $\tau$-tilting are red in color.

The number of support $\tau$-tilting $A$ modules is twenty four(see [1, Example 1.19]). Note that there are twelve support $\tau$-tilting $A$ modules in the middle line of the table and we cannot obtain the others by the construction.

{\bf Acknowledgments}\quad The author is grateful to Professor Bernhard Keller for answering questions about``the intermediate extension functor" by email. The author also thanks Professor Dong Yang and Professor Jiaqun Wei for their useful discussions. The author is indebted to the referees for helpful comments.

\section*{References}{ \footnotesize \baselineskip 13pt

\REF{[1]} Asai, S., Semibricks, {\it Int. Math. Res. Not}, 16(2020), 4993-5054.

\REF{[2]} Adachi, T., Iyama, O. and Reiten, I., $\tau$-tilting theory, {\it Compos. Math}, 150(2014), 415-452.

\REF{[3]} Angeleri H¨¹gel, L., Koenig, S. and Liu, Q., Recollements and tilting objects, {\it J. Pure Appl. Algebra}, 215(2011), 420-438.

\REF{[4]} Assem, I., Simson, D. and Skowro\'nski, A., Elements of the representation theory of associative algebras, Vol. 65(Cambridge University Press, Cambridge, 2006).

\REF{[5]} Bondarko, M.V., Weight structures vs. t-structures; weight filtrations, spectral sequences, and complexes(for motives and in general), {\it J.K-Theory}, 6(2010), 387-504.

\REF{[6]} Be\u\i linson, A.A., Bernstein, J. and Deligne, P.,
Faisceaux pervers, in Analysis and topology on singular spaces, I
(Luminy, 1981), 5--171, Ast\'erisque, 100, Soc. Math. France, Paris,
1982.

\REF{[7]} Chen, Q.H. and Lin, Y.N., Recollements of extension algebras, {\it Sci China Ser A}, 46(2003), 530-537.

\REF{[8]} Crawley-Boevey, W. and Sauter, J., On quiver Grassmannians and orbit closures for representation-finite algebras, {\it Math. Z}, 285(2017), 367-395.

\REF{[9]} Demonet, L., Iyama, O. and Jasso, G., $\tau$-tilting finite algebras, bricks and g-vector, {\it Int. Math. Res. Not}, 3(2019), 852-892.

\REF{[10]} Franjou, V. and Pirashvili, T., Comparison of abelian categories recollements, {\it Doc. Math}, 9(2004), 41-56.

\REF{[11]} Gabriel, P., Des cat$\acute{e}$gories ab$\acute{e}$liennes(French), {\it Bull. Soc. Math. France}, 90(1962), 323-448.

\REF{[12]} Gabriel, P., Indecomsable representations, II, Symposia Mathematica, Vol. ¢û(Convegno di Algebra Commutative, INDAM, Rome, 1971), pp. 81-104, Academic Press, London, 1973.

\REF{[13]} King, A., Moduli of representations of finite-dimensional algebras, {\it Q. J. Math}, 45(1994), 515-530.

\REF{[14]} Keller, B. and Scherotzke, S., Graded quiver varieties and derived categories, arXiv:1303.2318.

\REF{[15]} Liu, Q., Vit\'oria, J. and Yang, D., Gluing silting objects, {\it Nagoya Math. J}, 216(2014), 117-151.

\REF{[16]} Plamondon, P.G., $\tau$-tilting finite gentle algebras are representation-finite, {\it Pacific J. Math}, 302(2019), 709-716.

\REF{[17]} Psaroudakis, C., Homological theory of recollements of abelian categories, {\it J. Algebra}, 398(2014), 63-110.

\REF{[18]} Parra, C.E. and Vit\'oria, J., Properties of abelian categories via recollements, {\it J. Pure Appl. Algebra}, 223(2019), 3941-3963.

\REF{[19]} Psaroudakis, C. and Vit\'oria, J., Recollements of module categories, {\it Appl Categor Struct}, 22(2014), 579-593.

\REF{[20]} Ringel, C.M., Representations of K-species and bimodules, {\it J. Algebra}, 41(1976), 269-302.

\REF{[21]} Saor\'in, M. and Zvonareva, A., Lifting of recollements and gluing of partial silting sets, {\it Proc. Roy. Soc. Edinburgh Sect. A: Mathematics}, 152(2022), 209-257.

\REF{[22]} Zhang, P., Gorenstein-projective modules and symmetric recollements, {\it J. Algebra}, 388(2013), 65-80.
}

\end{sloppypar}

\end{document}